\documentclass[12pt,leqno,twoside,a4paper]{article}
\usepackage{amsfonts,amsmath,amsthm,amssymb}
\usepackage[english]{babel}
\newtheorem{theorem}{Theorem}[section]
\newtheorem{prop}[theorem]{Proposition}
\newtheorem{lemma}[theorem]{Lemma}
\newtheorem{corollary}[theorem]{Corollary}
\theoremstyle{definition}
\newtheorem{definition}[theorem]{Definition}
\theoremstyle{definitions}
\newtheorem{definitions}[theorem]{Definitions}
\newtheorem{example}[theorem]{Example}
\theoremstyle{remark}
\newtheorem{remark}[theorem]{Remark}
\theoremstyle{remarks}
\newtheorem{remarks}[theorem]{Remarks}
\numberwithin{equation}{section}

\newcommand{\de}{\delta}
\newcommand{\De}{\Delta}

\newcommand{\fl}{\longrightarrow}
\newcommand{\Fl}{\Longrightarrow}

\begin{document}
\title{\bf Noncommutative Symmetric functions \\ and $W$-polynomials}
\author{ {\bf  Jonathan
Delenclos and Andr\'{e} Leroy }\vspace{12pt}
\\ Universit\'{e} d'Artois,  Facult\'{e} Jean Perrin\\
Rue Jean Souvraz  62 307 Lens, France\\
     E-mail: delenclos@euler.univ-artois.fr\\
     E-mail: leroy@euler.univ-artois.fr}
\maketitle\markboth{ \bf J. DELENCLOS AND A. LEROY }{ \bf
Symmetric functions and $W$-polynomials}
%

\begin{abstract}
Let $K,S,D$ be a division ring an endomorphism and a
$S$-derivation of $K$, respectively.  In this setting we introduce
generalized noncommutative symmetric functions and obtain
Vi\`{e}te formula and decompositions of differential operators.
$W$-polynomials show up naturally, their connections with
$P$-independency, Vandermonde and Wronskian matrices are briefly
studied.  The different linear factorizations of $W$-polynomials
are analysed. Connections between the existence of LLCM of monic
linear polynomials with coefficients in a ring and the left duo
property are established at the end of the paper.
\end{abstract}

\section{Introduction}

Let $K$ be a division ring with center $k$.  Wedderburn
(Cf.\cite{We}) proved, among other results, that if $f(t)\in k[t]$
is irreducible but has a root $d_1$ in $K$, then $f(t)$ splits
linearly in $K[t]$ :$f(t)=(t-d_n)\cdots (t-d_1)$.  Moreover,
Wedderburn showed that if $f_i(t)$ stands for $(t-d_i)\cdots
(t-d_1)$ and $d$ is a conjugate of $d_1$ such that $f_i(d)\ne 0$,
then one can choose $d_{i+1}:=f_i(d)df_i(d)^{-1}$;.  In
particular, the elements $d_1,\dots,d_n$ are all conjugate to
$d_1$.  This result was extended by Jacobson using module theory
(Cf. \cite{J}) and extended to polynomials in an Ore extension
over a division ring by Lam and the second author.  Rowen and
Haile (Cf. \cite{HR}) studied in details factorizations of central
polynomials which are minimal polynomials of some elements in a
division rings.  They introduced good and very good factorizations
and gave applications to the structure of division rings.  Haile
and Knus used Wedderburn mathod to study algebras of degree 3 with
involutions 5\cite{HK}.  Rowen and Segev (\cite{RS}) also applied
the Wedderburn method for the study of the multiplicative group of
a finitedimensional division ring.

About 12 years ago, Gelfand and Retakh introduced noncommutative
symmetric functions via linear factorizations of quasideterminant
of Vandermonde matrices over division rings.  The form of these
symmetric functions already appears in the coefficients of the
above factorization of a central polynomial $f(t)=(t-d_n)\cdots
(t-d_1)$ if the elements $d_1,\dots,d_n$ are expressed in terms of
the roots of $f(t)$.  Gelfand and Retakh proved the required
symmetry property  and obtained Vi\`{e}te formula, Bezout and
Miura decompositions through computations with quasideterminants
(Cf \cite{GGRW}, \cite{GR}, \cite{GRW}). Later Wilson gave a
remarkable generalization of the fundamental theorem on symmetric
functions (Cf. \cite{W}, see also example \ref{free algebra
example}).

\noindent One of our aim in this paper is to show the link between
the Gelfand, Retakh, Wilson symmetric functions and the Wedderburn
polynomials (abbreviated $W$-polynomials) studied by Lam, Ozturk
and the second author.  The generalized Vandermonde matrices and
Wronskian matrices naturally show up in this context and enable us
to simplify presentation and generalize some of the results
mentioned above.  This is the content of section $2$ and $3$.

\noindent In section $4$, we show how to quickly introduce
quasideterminants in the different factorizations obtaining in
this way exactly the same expressions as the one given by Gelfand
and Retakh. We also compute the $LU$ decompositions of generalized
Vandermonde and Wronskian matrices.

\noindent The notion of $P$-independence was introduced in earlier
works by Lam and the second author (Cf.\cite{LL2} \cite{LL3} or
\cite{LL4}). This notion is recalled and studied briefly in
section $5$ in order to present the relations between Vandermonde
and Wronskian matrices and $W$-polynomials.  This leads to an
algorithmic test for $P$-independence (Cf. Proposition
\ref{algorithm for computing U}).  In fact, the $P$-independence
appears also in a somewhat hidden form in \cite{GRW}.  The results
of section $5$ give a new insight on the hypotheses used in this
latter reference (see remark \ref{remark on P-independence} b)).

\noindent In section $6$, we study all the factorizations of a
$W$-polynomial. These are shown to be in bijection with complete
flags in some spaces (Cf. Theorem \ref{bijection flags in E and
factorizations when all rooots in 1 CC } and Theorem
\ref{bijection factorizations and flags general case})

\noindent The last section is a first step towards considering
more general base ring than division rings.  It shows that a
natural choice could be to investigate the Wedderburn polynomials
and the symmetric functions with coefficients in a left duo or
more generally in an $(S,D)$ left duo ring.

\section{Noncommutative Symmetric functions }

\vspace{3mm}

Let $K,S,D$ be a division ring, an endomorphism of $K$ and a
$S$-derivation of $K$ respectively.  Throughout the paper we
denote by $R=K[t;S,D]$ the skew polynomial ring whose elements are
polynomials of the form $\sum_{i=0}^na_it^i$, where $a_0,\dots,a_n
\in K$.  The additive structure of $R$ is the usual one and
multiplication is bilinear and based on the commutation rule :
$$
ta=S(a)t + D(a) \quad \mbox {\rm{for} } a\in K.
$$
\noindent If $\{x_1,\dots,x_n\}$ are elements of $K$ we can
compute the least left common multiple $p_n(t)$ of the polynomials
$t-x_i$ in $R:=K[t;S,D]$. There are several ways of conducting the
computations and they will lead to different factorizations of the
polynomial $p_n(t)$ . We will denote $[f,g]_l$ (or simply $[f,g]$)
the monic polynomial which is a least left common multiple of $f$
and $g$ i.e. $[f,g]_l$ is a monic polynomial such that $Rf \cap Rg
= R[f,g]_l$.

For $f\in R$ and $a\in K$, we denote by $f(a)\in K$ the remainder
of $f$ right divided by $t-a$, i.e. $f(a)$ is the unique element
in $K$ such that $f-f(a)\in R(t-a)$.  This notion was introduced
in \cite{L'}.  For the sake of completeness let us recall the
important product formula.  In the sequel, for $a\in K$ and $c\in
K\setminus\{0\}$, we write $a^c$ for the $(S,D)$ conjugation of
$a$ by $c$ : $a^c:=S(c)ac^{-1}+D(c)c^{-1}$.
\begin{lemma}\rm {(Product formula)}
\label{product formula} Let $f, \, g \in R:= K[t;S,D]$ and $a \in
K$ :
\begin{enumerate}
\item $$ fg(a)= \left \lbrace
                 \begin{array}{l}
                 0 \quad \quad \quad \quad \quad \quad \; if \quad  g(a)=0 \\
                 f(a^{g(a)})g(a) \quad \quad if \quad g(a)\ne 0
                 \end{array}\right.
     $$
\item $$ [f,t-a]_l= \left \lbrace
                 \begin{array}{l}
                 \quad f \quad \quad \; \quad \quad \quad if \quad f(a)=0 \\
                 (t-a^{f(a)})f \quad \quad if \quad f(a)\ne 0
                 \end{array}\right. \quad\quad
     $$
\end{enumerate}
\end{lemma}
\begin{proof}

\noindent 1.  Right dividing $g(t)$ by $t-a$ we have
$g(t)=p(t)(t-a)+g(a)$ for some polynomial $p(t)\in R$.  Assuming
$g(a) \ne 0$, we have $f(t)= q(t)(t-a^{g(a)})+f(a^{g(a)})$ for
some $q(t)\in R$.  This gives $f(t)g(t)=
(q(t)(t-a^{g(a)})+f(a^{g(a)}))p(t)(t-a) +
q(t)(t-a^{g(a)})g(a)+f(a^{g(a)})g(a))$.  The fact that
$(t-a^{g(a)} )g(a)=S(g(a))(t-a)$ then leads quickly to the
conclusion.

\noindent 2. Is is easy to check that, if $f(a)\ne 0$,
$(t-a^{f(a)})f\in R(t-a)$.  This gives the result.
\end{proof}

Let us notice that $deg([f,t-a]_l)\le deg(f)+1$ and that equality
occurs if and only if $f(a)\ne 0$.  This will be used freely in
the paper.

 \noindent For $\{x_1,\dots,x_n\}\subseteq K$, we can now construct
 the $llcm$ of $t-x_1,\dots ,t-x_n \in R=K[t;S,D]$.

\begin{example}
\label{example quadratic W-polynomial} {\rm Suppose $x_1\ne x_2$
are elements in $K$. We have :
$$[t-x_1,t-x_2]_l = (t-x_1^{x_1 - x_2})(t-x_2)=(t-x_2^{x_2 -
x_1})(t-x_1)$$

Comparing coefficients of degree $0$ and $1$, this immediately
leads to the following generalized symmetric functions :

$$\Lambda_1(x_1,x_2) = x_1^{x_1 - x_2}+S(x_2) =x_2^{x_2 - x_1} + S(x_1)$$

$$\Lambda_2(x_1,x_2) = x_1^{x_1 - x_2}.x_2-D(x_2)=x_2^{x_2 - x_1}.x_1-D(x_1) $$
}
\end{example}

\vspace{5mm}

\noindent In order to exhibit the symmetric functions in
$\{x_1,\dots,x_n\}\subseteq K$ with $n\ge 2$, we introduce some
notations and a definition:

\noindent Let us put $p_j=[t-x_i \, |\, i\le j]_l$ for
$j=1,\dots,n$. It is useful to also define $p_0:=1$. Using the
Lemma \ref{product formula} $2.$, we then get

$$ p_n(t)= \left \lbrace
                 \begin{array}{l}
                 p_{n-1}(t) \quad if \quad p_{n-1}(x_{n})=0 \\
                 (t-x_n^{p_{n-1}(x_{n})})p_{n-1}(t) \quad if \quad p_{n-1}(x_{n})\ne 0
                 \end{array}\right.
                 $$

This enables us to compute $p_n$ by induction on $n$ starting with
$p_0(t)=1$.

Notice that $x_1,\dots,x_n$ are right roots of $p_n(t)$ and the
above formula means that $p_i(t)$ can be computed from
$p_{i-1}(t)$ by requiring that $x_i$ is a root of $p_i(t)$. It is
obvious that $p_n(t)$ can be computed in different ways depending
on the order in which the roots are added.  In example
\ref{example quadratic W-polynomial} we have seen the $2$
different expressions for $p_2(t)=[t-x_1,t-x_2]$ when $x_1\ne
x_2$.
\begin{definition}
\label{P-independency} We say that the set
$\{x_1,\dots,x_n\}\subseteq K$ is $P$-independent if $deg([t-x_i |
1\le i \le n ])=n$.
\end{definition}

We leave the easy proof of the next lemma to the reader.

\begin{lemma}
\label{first equivalence for P independence}
$\{x_1,\dots,x_n\}\subseteq K$ is $P$-independent if and only if
for any $i\in \{1,\dots,n-1\}$, $p_i(x_{i+1})\ne 0$.
\end{lemma}

\noindent   An equivalent definition was given in \cite{LL2}. This
notion will be studied more deeply in section $5$ (Cf. Theorem
\ref{characterization of W-polynomials}). Of course, if
$\{x_1,\dots,x_n\}$ is a $P$-independent set then for every $i\le
n $ the set $\{x_1,\dots,x_i\}$ is $P$-independent as well. In
this case we put
$$y_i=x_i^{p_{i-1}(x_i)} \; {\mbox{\rm for} }\;
i\in\{1,\dots,n\}$$
We then have:
$$
 [t-x_j \,|\,
j=1,\dots,i]_l=p_i(t)=(t-y_i)\dots(t-y_1) \; {\mbox{\rm for} }\;
i\in\{1,\dots,n\} \quad(A)
$$

\noindent We define $(-1)^k\Lambda_k^i(X), \, 0 \le k \le i$, to
be the coefficient of degree $i-k$ of $p_i(t)$ i.e.
$p_i(t)=\sum_{k=0}^i(-1)^k\Lambda_k^it^{i-k}$.  As in example
\ref{example quadratic W-polynomial}, introducing the roots of
$p_i$ in different orders give different factorizations of this
polynomial and this shows that the coefficients $\Lambda_k^i(X)$
are symmetric with respect to permutations of $X$. In terms of the
$y_j$'s they can be written as :

\noindent$ \Lambda_0^1=1.\\
\Lambda_1^1 = y_1. \\
\Lambda_0^2 = 1. \\
\Lambda_1^2 = y_2 + S(y_1). \\
\Lambda_2^2=y_2y_1 -D(y_1).\\$

Assume that the $\Lambda_0^i,\dots,\Lambda_i^i$ have been defined
(with $i<n$).  Remarking that $p_{i+1}(t)=(t-y_{i+1})p_i(t)$ we
have:

\noindent $\Lambda_0^{i+1}=1.\\
\Lambda_1^{i+1}= y_{i+1} + S(\Lambda_1^i). \\
\Lambda_2^{i+1}=y_{i+1}\Lambda_1^i +S(\Lambda_2^i)-D(\Lambda_1^i).\\
\hspace{3mm} \vdots \quad \quad \quad \quad \vdots\\
\Lambda_{k}^{i+1} = y_{i+1}\Lambda_{k-1}^i +
S(\Lambda_k^i)-D(\Lambda_{k-1}^i).\\
\vdots \quad\quad\quad\quad \vdots \\
\Lambda_{i+1}^{i+1}= y_{i+1}\Lambda_i^i -D(\Lambda_{i-1}^i).$

\vspace{2mm}
\begin{remark}
\label{the classical symmetric functions}
\begin{enumerate}
\item[a)] In the classical case ($S=Id.,\, D=0$) one can easily
describe the $\Lambda$'s in terms of the $y_i$'s.  From (A) above
one gets :

\noindent$\Lambda_0^n=1.\\
\Lambda_1^n(X) = y_1+\dots y_n. \\
\Lambda_2^n(X) =\sum_{i<j}y_jy_i. \\
\Lambda_3^n(X)=\sum_{i<j<k}y_ky_jy_i.\\
\vdots \quad \quad \quad \quad \quad\vdots \\
\Lambda_n^n(X)=y_ny_{n-1}\dots y_1.$.
\item[b)] Let us notice that, if $x_1,\dots,x_n\in K$, an algorithm can be
written to check $P$-independency of these elements and in this
case to compute the elements $y_1,\dots,y_n$ and the
$\Lambda^i_k$.
\item[c)]In their works on noncommutative symmetric functions
Gelfand, Rethak and Wilson used quasideterminants to obtain (in
case when $S=Id.$ and $D=0$ ) the above symmetric functions.  With
this point of view the symmetry of these functions was not easy to
prove and somewhat surprising  (see comments in the introduction
of \cite{GRW}).  With the above approach the symmetry is clear
even in the $(S,D)$-setting.
\end{enumerate}
\end{remark}

Let us end this section with an important example and result
extracted from \cite{W}.

\begin{example}
\label{free algebra example} {\rm Let $A:=k<x_1,\dots,x_n>$ be a
free algebra over a commutative field $k$ in $n$ noncommutative
free variables. Let us denote by $K$ the universal field of
fractions of $A$ (aka. the free field in $n$ noncommutative
variables) and consider the (usual) polynomial ring $K[t]$. The
least left common multiple $p(t):=[t-x_1,\dots,t-x_n]_l$
factorizes in $K[t]$ as above : $p(t)=(t-y_n)\cdots(t-y_n)$.  In
this context an analogue of the classical fundamental theorem on
symmetric functions was given in \cite{W} : If a polynomial $f\in
k[y_1,\dots,y_n]$ is symmetric with respect to  the variables
$x_1,\dots,x_n$, then $f\in
k[\Lambda_1^n(X),\dots,\Lambda_n^n(X)]$, where for $1\le i \le n$,
the elements $\Lambda_i^n(X)$ are described in Remark \ref{the
classical symmetric functions} a).

 }
\end{example}

\vspace{5mm}
\section{Vi\`{e}te, Bezout and Miura  decompositions}

\vspace{3mm}

Suppose $\{x_1,x_2,\dots,x_n\}\subseteq K$ is a $P$-independent
set.  With the above notations, we have :

$$
[t-x_i \,|\,
i=1,\dots,n]_l=p_n(t)=\sum_{i=0}^n(-1)^i\Lambda_i^nt^{n-i}
\quad\quad (B)
$$

This is the Vi\`{e}te formula expressing the coefficients of a
polynomial in terms of its roots, assumed to be $P$-independent.
 This condition will be studied in section $5$.
In case when $S=Id.$ and $D=0$ this expression was obtained by
Gelfand and Retakh (Cf. \cite{GGRW}) using quasideterminants
techniques.

If $z\in K$ is not a root of $p_n$, we get from equation (A) in
section $2$,
$p_n(z)=((t-y_n)p_{n-1})(z)=(z^{p_{n-1}(z)}-y_n)p_{n-1}(z)$.
Putting $z_i=z^{p_{i-1}(z)}$ for $1\le i \le n$ ($p_0:=1$), an
easy induction leads to:
$$
p_n(z)=(z_n-y_n)(z_{n-1}-y_{n-1})\dots (z_1-y_1) \quad\quad (C)
$$

This is sometimes called the Bezout decomposition (Cf.
\cite{GGRW}).

\noindent For $a\in K$, the left $R$-module $R/R(t-a)$ induces a
left $R$-module structure on $K$.  In particular, since for $x\in
K $, we have $tx=S(x)t+D(x)=S(x)(t-a)+S(x)a+D(x)$ and we conclude
that the action of $t\in R$ on $x\in K$ is given by
$t.x=S(x)a+D(x)$.  The action of $t$ on $K$ will be denoted by
$T_a$; we thus have $T_a(x)=S(x)a+D(x)$.  In general the structure
of $_RK$ is given by computing remainders of right division by
$t-a$.  We thus have for $f(t)\in R$ and $x\in K$,
$f(t).x=f(T_a)(x)$. For easy reference we give the following
explicit definitions :

\begin{definitions}
\label{PLT, Conju. class, C(a)} \hfill
\begin{enumerate}
\item[a)] $T_a: K\fl K : x\mapsto S(x)a + D(x)$.
\item[b)] For $a\in K$ we define the $(S,D)$-conjugacy class of $a$ by $\Delta(a):=\{a^x\,|\, x\in K
\setminus\{0\}\}$.
\item[c)] For $a\in K$, $C^{S,D}(a):=\{x\in K\setminus\{0\} \,|\, a^x=a\}\cup
\{0\}$ stands for the $(S,D)$-centralizer of $a$.
\item[d)] For $a\in K$ and $i\ge 0$, we put $N_i(a)=t^i(a)$.
\end{enumerate}
\end{definitions}

In the next lemma we collect some properties of the map $T_a$.

\begin{lemma}
\label{properties of T_a} Let $f(t)=\sum_{i=0}^na_it^i$ be any
element of $R=K[t;S,D]$ and $a\in K$. Then
\begin{enumerate}
\item $f(a)=\sum_{i=0}^na_iN_i(a)$.
\item $f(T_a)(x)=f(a^x)x$.
\item $C^{S,D}(a)$ is a subdivision ring of $K$.
\item $f(T_a)$ is a right $C:=C^{S,D}(a)$-linear map i.e. $f(T_a)\in Hom_C(K_C)$.
\item $Ker(f(T_a))=\{y\in K\setminus\{0\} \,|\, f(a^y)=0\}\cup \{0\}$ is a right
$C:=C^{S,D}(a)$ vector space denoted $E(f,a)$.
\end{enumerate}
\end{lemma}
\begin{proof}
\noindent 1.  $f(a)=\sum_{i=0}^na_it^i(a)=\sum_{i=0}^na_iN_i(a)$.

\noindent 2.  By the paragraph preceding the definitions in
\ref{PLT, Conju. class, C(a)} we have, for $f(t)\in R$ and $x\in
K$, $f(t).x=f(T_a)(x)$ where the left hand side refers to the
action of $f(t)$ on $x$ considered as an element of $R/R(t-a)$;
i.e. $f(t).x$ stands for the remainder of $f(t)x$ after right
division by $t-a$. In other words, $f(t).x=(f(t)x)(a)=f(a^x)x$.

\noindent 3. This is easy to check.

\noindent 4. We also leave this to the reader.

\noindent 5. It is clear from 4. above that $E(f,a)$ is a right
$C$-vector space.  Later (Cf.\ref{dim E(f) and deg(f)}) we will
show that $dim_C(E(f,a))\le deg (f)$.
\end{proof}

The map $T_a$ is in fact an $(S,D)$-pseudo-linear transformation.
For more information on these maps we refer to (\cite{L}). Notice
also that in earlier papers the map $f(T_a)$ above was denoted
$\lambda_{f,a}$ (Cf. e.g. \cite{LL3}). Remark also that the last
statement in \ref{properties of T_a} says that the nonzero
solutions of the operator equation $f(T_a)=0$ are the roots of
$f(t)$ belonging to $\Delta (a)$.  It is thus natural to look at
the special case when all the roots belong to one singular
$(S,D)$-class $\Delta(a):=\{a^x \,|\,x\in K\setminus\{0\} \}$.  In
particular, $\Delta (0)=\{D(y)y^{-1}\,|\, y\in K\setminus\{0\}\}$
is the set of so called logarithmic derivatives and the nonzero
solutions of a differential polynomial $f=D^n+a_{n-1}D^{n-1}+\dots
+a_1D+a_0$ coincide with the roots of $f(t)$ belonging to $\Delta
(0)$ (Cf. \cite{L'}).

Assume that $\{x_1=a^{u_1},\dots,x_n=a^{u_n}\}$ are
$P$-independent we then have :

\noindent
$y_i=x_i^{p_{i-1}(x_i)}=(a^{u_i})^{p_{i-1}(a^{u_i})}=a^{p_{i-1}(a^{u_i})u_i}=a^{w_i}$
where, as above, $p_{i-1}=[t-a^{u_j}\,|\, 1\le j \le i-1 ]_l$ and
$$w_i=p_{i-1}(a^{u_i})u_i=p_{i-1}(T_a)(u_i)$$.

Equation (A) in section $2$ takes now the form

$$
p_n(t)= (t-a^{w_n})\dots (t-a^{w_1}) \quad\quad \quad \quad (D)
$$

Applying this decomposition to the pseudo-linear transformation
$T_a$ we get :

$$
p_n(T_a)=(T_a-a^{w_n})\dots (T_a-a^{w_1}) \quad\quad \quad (E)
$$

This is the Miura decomposition.  This was obtained in the case
when $a=0$ (i.e. for $T_0=D$) and $S=Id.$ in \cite{GRW}, section
$4$ using quasideterminants techniques.

Although the equations we obtained are independent of the
quasideterminants it is worth to show how we can make them appear.
This is one of the objectives of the next section.

\section{Quasideterminants, (S,D)-Vandermonde and
Wronskian matrices}

For a matrix $A=(a_{ij})\in M_n(R)$, where $R$ is a ring, there
might be up to $n^2$ quasideterminants denoted by $|A|_{ij}$.
$|A|_{ij}$ is defined when the matrix $A^{ij}$, obtained from $A$
by deleting the $ith$ row and the $jth$ column, is invertible. In
this case we have:

$$
|A|_{ij} = a_{ij}-r^i_j(A^{ij})^{-1}c_i^j\,,
$$
where $r^i_j$ is the $ith$ row of $A$ from which the $jth$ element
has been suppressed and the $c_i^j$ is the $jth$ column of $A$
from which the $ith$ element has been suppressed.

We need the following special case of a more general result :
\begin{lemma}
\label{QD_{nn} } If $A\in M_n(K)$ is invertible and $A^{nn}$ is
invertible then
$$|A|_{nn}=((A^{-1})_{nn})^{-1}.$$
\end{lemma}
\begin{proof}
Let us write $$A=\begin{pmatrix}
                 A^{nn}        & c    \\
                 l        & a_{nn}    \\
   \end{pmatrix} \quad {\mbox {\rm and} } \quad
   A^{-1}=\begin{pmatrix}
                 B        & v     \\
                 u       & \alpha  \\
          \end{pmatrix}\,,
$$
where $c,v$ are columns, $l,u$ are lines and $a_{nn},\alpha \in
K$. Comparing the last column on both side of $AA^{-1}=Id.$ we get
two equations : $A^{nn}v+c\alpha=0$ and $l.v +a_{nn}\alpha=1$.
These lead to $(a_{nn}-l(A^{nn})^{-1}c)\alpha=1$ which shows that
$\alpha^{-1}=|A|_{nn}$.
\end{proof}

Since we are working in an $(S,D)$-setting we introduce a
generalized Vandermonde matrix as follows (Cf. \cite{LL1}): let
$\{x_1,\dots,x_n\}$ be any subset of $K$,

$$
V_n^{S,D}(x_1,\dots,x_n)=\begin{pmatrix}
                 1          & 1      & \dots  & 1  \\
                 x_1        & x_2    & \dots  & x_n \\
                 N_2(x_1)   &N_2(x_2)& \dots  & N_2(x_n)\\
                 \dots      & \dots  & \dots  & \dots \\
                 N_{n-1}(x_1)&N_{n-1}(x_2)& \dots &N_{n-1}(x_n)
                 \end{pmatrix}.
$$

Recall that, for $x\in K$, $N_i(x)$ is the evaluation $t^i(x)$. Of
course, the $N_i$'s can be computed independently of the
evaluation process i.e. in terms of $S$ and $D$.  Indeed, using
Lemma \ref{product formula} for $t^{i+1}=tt^i$, one has
$N_{i+1}(x)= x^{t^i(x)}t^i(x)=x^{N_i(x)}N_i(x)=S(N_i(x))x +
D(N_i(x))$. This gives a recurrence formula for the computation of
$N_i(x)$. In the classical setting ($S=Id.$ and $D=0$) one has
$N_i(x)=x^i$.  The evaluation of $f(t)=\sum_{i=0}^{n-1}a_it^i\in
K[t;S,D]$ at $x\in K$ can now be expressed as
$f(x)=\sum_{i=0}^{n-1}a_iN_i(x)$.  Hence, as in the classical case
$(S=Id.,D=0)$, we have, for $\{x_1,\dots,x_n\}\subseteq K$,

$$
(f(x_1),f(x_2),\dots ,
f(x_n))=(a_0,a_1,\dots,a_{n-1})V(x_1,\dots,x_n).
$$

In \cite{LL1} we also computed the inverse of an (invertible !)
generalized Vandermonde matrix.  Let us recall this :

\noindent for $1\le i \le n$, define $g_i(t):=[t-x_j \, | \, 1\le
j \le n , \, j\ne i ]_l$.  Since, for every $1\le i \le n$, $deg
(g_i(t) )\le n-1$ one can write
$g_i(t)=\sum_{j=0}^{n-1}c_{ij}t^j$. Consider the matrix
$C=(c_{ij})$. Since $g_i(x_j)=0$ if $i\ne j$, one has
$$
CV_n^{S,D}(x_1,\dots,x_n)=diag(g_1(x_1),\dots,g_n(x_n)) \quad
\quad (F)
$$

\begin{prop}
\label{QD of VDM}
 Suppose that the elements $x_1,\dots,x_n$ are
$P$-independent. Then the Vandermonde matrix
$V=V_n^{S,D}(x_1,\dots,x_n)$ is invertible and $|V|_{nn}=
g_n(x_n)$, where $g_n=[t-x_i \, | \, i=1,\dots,n-1]_l$.
\end{prop}
\begin{proof}
Let us recall,from section $2$, that $p_n(t):=[t-x_i \,|_l\,
i=1,\dots,n]$.  Assuming that the elements $x_1,\dots,x_n$ are
$P$-independent, we have that $n=deg(p_n(t))=deg [t-x_i,g_i(t)]$,
for $1\le i\le n$.  This implies that $deg(g_i(t) )= n-1$ and
$g_i(x_i)\ne 0$.  Hence,

$$(V(x_1,\dots,x_n))^{-1} =
diag(g_1(x_1)^{-1},\dots,g_n(x_n)^{-1})C .$$

Since $g_n(t)$ is a monic polynomial we have $C_{n,n-1}=1$ and so
$(V^{-1})_{n,n}=g(x_n)^{-1}$. Lemma \ref{QD_{nn} } shows that
$|V|_{nn}=g_n(x_n)$.
\end{proof}

Let us now consider the case of generalized Wronskian matrices.
Let us recall that, for $x\in K,\; T_a(x)=S(x)a+D(x)$ (Cf.
\ref{PLT, Conju. class, C(a)}). For $u_1,\dots,u_n \in K$, we put
:
$$
W_{n,a}^{S,D}(u_1,\dots,u_n)=\begin{pmatrix}
                 u_1        & \dots    & u_n  \\
                 T_a(u_1)   & \dots   & T_a(u_n) \\
                 T_a^2(u_1) & \dots   & T_a^2(u_n)\\
                 \dots      & \dots   & \dots \\
                 T_a^{n-1}(u_1)&\dots &T_a^{n-1}(u_n)
                 \end{pmatrix}.
$$

\begin{example}
{\rm It is worth to mention the special form of the Wronskian
matrix when $a=0,\;S=Id.$.  The corresponding Wronskian matrix is
of the form :
$$
W_{n,0}^{Id.,D}(u_1,\dots,u_n)=\begin{pmatrix}
                 u_1        & \dots    & u_n  \\
                 D(u_1)   & \dots   & D(u_n) \\
                 D^2(u_1) & \dots   & D^2(u_n)\\
                 \dots      & \dots   & \dots \\
                 D^{n-1}(u_1)&\dots &D^{n-1}(u_n)
                 \end{pmatrix},
$$
and when $a=1,\;D=0$.  The corresponding Wronskian matrix is of
the form :

$$
W_{n,1}^{S,0}(u_1,\dots,u_n)=\begin{pmatrix}
                               u_1        & \dots    & u_n  \\
                               S(u_1)   & \dots   & S(u_n) \\
                               S^2(u_1) & \dots   & S^2(u_n)\\
                               \dots      & \dots   & \dots \\
                                S^{n-1}(u_1)&\dots &S^{n-1}(u_n)
                               \end{pmatrix}.
$$
}
\end{example}

The next lemma establishes connections between the Vandermonde and
Wronskian matrices and compute a quasideterminant of this
Wronskian.

\begin{prop}
\label{QD of WRSK} Let $u_1,\dots,u_n$ be nonzero elements in $K$.
For any $a\in K$ one has
\begin{enumerate}
\item
$V_n^{S,D}(a^{u_1},\dots,a^{u_n})diag(u_1,\dots,u_n)=W_{n,a}^{S,D}(u_1,\dots,u_n)$.
\item $V_n^{S,D}(a^{u_1},\dots,a^{u_n})$ is invertible if and only
if $W_{n,a}^{S,D}(u_1,\dots,u_n)$ is invertible.
\item If $W:=W_{n,a}^{S,D}(u_1,\dots,u_n)$ is invertible then we
have
$$|W|_{n,n}=g_n(a^{u_n})u_n=g_n(T_a)(u_n)$$
where $g_n=[t-a^{u_1},\dots,t-a^{u_{n-1} }]_l$.
\end{enumerate}
\end{prop}
\begin{proof}
\noindent $1)$   Recall that for $f(t)\in K[t;S,D]$ we have proved
in \ref{properties of T_a} that $f(T_a)(x)=f(a^x)x$.  In
particular, $N_i(a^x)x=t^i(a^x)x=T_a^i(x)$.  From this, one gets
easily the equation relating Vandermonde and Wronskian matrices.

\noindent $2)$ This is an obvious consequence of $1)$

\noindent $3)$ From $1)$ we obviously get that
$(W^{-1})_{nn}=u_n^{-1}(V_n^{S,D}(a^{u_1},\dots,a^{u_n})^{-1})_{nn}$.
Lemma \ref{QD_{nn} } yields the result.
\end{proof}
 These propositions show that the form of Viì\`{e}te,
Bezout and Miura decompositions obtained in \cite{GGRW},
\cite{GR}, \cite{GRW} are special cases (S=Id. D=0) of the one
obtained in section $3$. Indeed, it suffices to replace the
evaluations of least left common multiples by the corresponding
quasideterminants to find back the formulas from the papers
mentioned above. To be more precise, let us write
$V_i:=|V_i^{S,D}(x_1,\dots,x_i)|_{i,i}$ then, thanks to
Proposition \ref{QD of VDM}, the $y_i$'s obtained in equation (A)
can be written $y_i=x_i^{V_i}$. Similarly writing
$Z_i:=|V_i^{S,D}(x_1,\dots,x_{i-1},z)|_{ii}$ the $z_i$'s appearing
in (C) are given by $z_i=z^{Z_i}$.  Finally, we remark that the
$w_i$'s in (D) are equal to
$w_i=|W_{i,a}^{S,D}(u_1,\dots,u_i)|_{i,i}$.
\begin{remark}
The fact that the $\Lambda_i^k$ are symmetric in $x_1,\dots,x_n$
is somewhat surprising when the $y_i$'s are expressed with
quasideterminants (i.e. $y_i=x_i^{V_i}$, as above).  This fact is
obvious from the point of view of least left common multiple as
developed in section $2$.

\end{remark}

We end this section with two easy propositions giving
$LU$-decomposition (this corresponds also to the strict Bruhat
normal form) of invertible Vandermonde and Wronskian matrices. Let
us first recall that a matrix $A\in GL_n(K)$ can be written in the
form $A=LDPU$ where $L$ is lower unitriangular $D$ is diagonal $P$
is a permutation matrix and $U$ is an upper unitriangular matrix.
(see \cite{D} p. $128$, for details). This form is unique and has
some importance nowadays due to its usage in computer packages for
solving linear systems of equations.  It will be easier for the
notations to index the rows of the matrices starting with $0$.
Notice that, in this case, the diagonal elements of a matrix
$a_{i,j}$ are $a_{0,1},\dots,a_{n-1,n}$.  Let us first state,
without proof, the following easy lemma :
\begin{lemma}
\label{AV}

Let $\{x_1,\dots,x_{n},a\}$ be a subset of $K$ and
$\{u_1,\dots,u_n\}$ be a subset of $K\setminus\{0\}$.  Let $l\ge
0$ be an integer and $A=(a_{ij})\in M_{(l+1)\times n}(K)$.  Then
$$
AV_n^{S,D}(x_1,\dots,x_n)=(f_i(x_j)) \quad {\mbox{\rm and } }
\quad AW_{n,a}^{S,D}(u_1,\dots,u_n)=(f_i(T_a)(u_j))
$$
where for $0\le i \le l,\; f_i(t):=\sum_{j=1}^{n}a_{ij}t^{j-1}$.
\end{lemma}

Let $\{x_1,\dots,x_n\}\subset K$ be a $P$-independent subset of
$K$. Recall from section $2$ that $p_i(t)=[t-x_j\,|\,j\le
i]_l=\sum_{k=0}^i(-1)^k\Lambda_k^it^{i-k}$.   We define
$\Lambda\in M_n(K)$ a lower unitriangular matrix via
$$
\Lambda=\begin{pmatrix}
                 1=\Lambda_0^0   &   0         &  0   &\dots    & 0 \\
                 -\Lambda_1^1  &1=\Lambda^1_0 & 0    &\dots    & 0 \\
                 \Lambda_2^2   &-\Lambda^2_1  & 1=\Lambda^2_0&\dots& 0\\
                 \dots      & \dots  & \dots  & \dots &\dots\\
                 (-1)^{n-1}\Lambda_{n-1}^{n-1}&(-1)^{n-2}\Lambda^{n-1}_{n-2}&
                 (-1)^{n-3}\Lambda^{n-1}_{n-3}&\dots& 1
                 \end{pmatrix}
$$

In the following proposition it is important to notice that the
rows of the matrix $U$ are indexed starting with $0$ i.e. the
first row of $U$ is $(U_{01},\dots,U_{0n})$.
\begin{prop}
\label{LU decomposition of V and W} Let $\{x_1,\dots,x_n\}\subset
K $ be a $P$-independent set.  With the above notations we have:
\begin{enumerate}

%

\item $\Lambda V=U$ where $U_{ij}=\left \lbrace
                 \begin{array}{l}
                \quad 0 \quad \quad if \quad i\ge j \\
                 p_i(x_j) \quad if \quad i< j
                 \end{array}
                 \right. $
\item $V=\Lambda^{-1}U=\Lambda^{-1}diag(1,p_1(x_2),\dots,p_{n-1}(x_n))U'$
where $U'$ is an upper unitriangular matrix.
\end{enumerate}
\end{prop}
\begin{proof}
\noindent 1.  This is obvious since the elements of the
$i^{th}$-row of the matrix $\Lambda$ are the coefficients of the
polynomial $p_i(t)=[t-x_j \,|\,j\le i]_l$ for $i=0,\dots,n-1$
($p_0=1$, see section $2$ or equation (B) in section $3$).  The
above Lemma \ref{AV} then yields the result.

\noindent 2. Obviously $\Lambda$ is invertible and the
$P$-independence of $x_1,\dots,x_n$ shows that $p_i(x_{i+1})\ne
0$. We can thus define
$U'=diag(1,p_1(x_2)^{-1},\dots,p_{n-1}(x_{n+1})^{-1})U$ and the
last equality follows.

\end{proof}

Since the Vandermonde and Wronskian matrices are so closely
related (Cf. \ref{QD of WRSK}) it is not surprising that we get a
similar result for Wronskian matrices.

\begin{prop}
Let $\{u_1,\dots,u_n\}$ be a subset of nonzero elements in $K$
which are right independent over $C^{S,D}(a)$ for some $a\in K$.
Then
\begin{enumerate}
\item $\Lambda W=Z$ where $Z_{ij}=\left \lbrace
                 \begin{array}{l}
                \quad 0\quad \quad \quad \quad if \quad i>j \\
                 p_i(T_a)(u_j) \;\quad if \quad i\le j
                 \end{array}
                 \right. \quad$
\item $W=(\Lambda)^{-1}Z=\Lambda^{-1}diag(1,p_1(T_a)(u_2),\dots,p_{n-1}(T_a)(u_n))Z'$
where $Z'$ is an upper unitriangular matrix.
\end{enumerate}
\end{prop}
\begin{proof}
The proof, based on the above lemma \ref{AV}, is similar to the
proof of the previous proposition.
\end{proof}

\begin{remark}
\label{LU in general U algorithmically}
\begin{enumerate}
\item[a)] In fact one can get an $LU$ decomposition of a
Vandermonde (resp. Wronskian) matrix without assuming that the set
$\{x_1,\dots,x_n\}$ is $P$-independent (resp. $\{u_1,\dots,u_n\}$
is right $C^{S,D}(a)$ linearly independent). Indeed there always
exist monic polynomials $q_i=\sum_jq_{ij}t^j$ of degree $i$ such
that $q_i(x_j)=0$ for $1\le j\le i$. The matrix $L =(q_{ij})$ of
the coefficients of these polynomials will give the invertible
lower unitriangular matrix $L$ and $U=LV$ will be an upper
triangular matrix.
\item[b)] The matrix $U$ can be algorithmically computed and
offers a way for testing the invertibility of a Vandermonde
matrix. This will be explained at the end of section $5$.
\end{enumerate}

\end{remark}

\section{P-independence and $W$-polynomials}

Let us start this section with the formal definition of a
Wedderburn polynomial.
\begin{definition}
A monic polynomial $f(t)\in R=K[t;S,D]$ is a Weddderburn
polynomial if there exists $A\subseteq K$ such that
$Rf(t)=\cap_{a\in A}R(t-a)$.
\end{definition}

The $W$-polynomials were introduced in \cite{LL3} and studied in
depth in \cite{LL4} and \cite{LLO}.  They are special kind of
fully reducible polynomials (also called completely reducible
polynomials by Ore ) Cf.\cite{C}, \cite{LO} for more details. From
the definition above, it is clear that the $W$-polynomials are
exactly the polynomials we have been dealing with from the
beginning of this paper.   In particular, these are exactly the
polynomials for which we have presented the symmetric functions an
the Vi\`{e}te formulas in section $2$ and $3$. The different
factorizations of these polynomials will be presented in later
sections. In the above mentioned papers the second author in
collaboration with Lam and Ozturk studied the roots, the
factorizations, the products of $W$ -polynomials.  It is also
clear that these polynomials are related to $(S,D)$-algebraic sets
and a lot of information contained in related works are relevant
to $W$-polynomials.   The interested reader can refer to the
bibliography mentioned in \cite{LL3} \cite{LL4}, or \cite{LO}. Let
us recall, without proofs, a few characterizations of these
polynomials :
\begin{prop}
\label{characterization of W-polynomials}

For a monic polynomial $f$ of degree $n$ in $R=K[t;S,D]$, the
following are equivalent :
\begin{enumerate}
\item[i)] $f$ is a $W$-polynomial.
\item[ii)] There exists a subset $\{a_1,\dots,a_n\}$ in $K$ such that
$Rf=\cap_{i=1}^n R(t-a_i)$.
\item[iii)] There exists a subset $\{a_1,\dots,a_n\}$ in $K$ such that
$V=V_n^{S,D}(a_1,\dots,a_n)$ is invertible and
$$
C_fV=S(V)diag(a_1,\dots,a_n)+D(V) .$$
\item[iv)] The companion matrix $C_f$ is $(S,D)$-diagonalisable i.e.
there exists an invertible matrix $P$ in $M_n(K) $ and a diagonal
matrix $\Delta$ such that $S(P)C_fP^{-1}+D(P)P^{-1}=\De$.

\end{enumerate}
\end{prop}

\noindent Obviously the $W$-polynomials are strongly related to
the notion of $P$-independency as defined in section $2$
definition \ref{P-independency}. We give now a few
characterizations of this notion. Let us recall, from definitions
\ref{PLT, Conju. class, C(a)}, that $C^{S,D}(a)=\{0\ne x \in K \,
| \, a^x=a\}\cup\{0\}$ is a subdivision ring of $K$ called the
$(S,D)$-centralizer of $a$.

\begin{theorem}
\label{equivalence for P-independence}

Let $\{x_1,\dots,x_n\}$ be a subset of $K$.  The following are
equivalent :
\begin{enumerate}
\item[i)] $\{x_1,\dots,x_n\}$ is a $P$-independent set.
\item[i')] For every subset $\{i_1,\dots,i_r\}\subseteq
\{1,\dots,n\}$, $\{x_{i_1},\dots,x_{i_r}\}$ is $P$-independent.
\item[ii)] $V(x_1,\dots,x_n)$ is invertible.
\item[ii')] For every subset $\{i_1,\dots,i_r\}\subseteq
\{1,\dots,n\},\;V(x_{i_1},\dots,x_{i_r})$ is invertible.
\end{enumerate}
\noindent If there exists $a\in K$ such that
$\{x_1,\dots,x_n\}\subseteq \Delta (a)$ say $x_i=a^{u_i}$ for
$i=1,\dots,n$ and $\{u_1,\dots,u_n\}\subseteq K\setminus \{0\}$.
Then the above statements are also equivalent to the following
ones :
\begin{enumerate}
\item[iii)] $\{u_1,\dots,u_n\}$ is right
$C^{S,D}(a)$-independent.
\item[iii')] For every subset $\{i_1,\dots,i_r\}\subseteq
\{1,\dots,n\},\;\{u_{i_1},\dots,u_{i_r})$ is right
$C^{S,D}(a)$-independent.
\item[iv)] $W_a^{S,D}(u_1,\dots,u_n)$ is invertible.
\item[iv')] For every subset $\{i_1,\dots,i_r\}\subseteq
\{1,\dots,n\},\;W_a^{S,D}(u_{i_1},\dots,u_{i_r})$ is invertible.
\end{enumerate}
\end{theorem}
\begin{proof}
\noindent  $i)\Fl ii)$. This was already proved in \ref{QD of
VDM}.

\noindent  $ii)\Fl i)$. If $V_n^{S,D}(x_1,\dots,x_n)$ is
invertible but $deg([t-x_i\,|\,i=1,\dots,n]_l)<n$ then we claim
that there exists $i\in \{1,\dots,n\}$ such that $g_i(x_i)=0$ (let
us recall that $g_i(t):=[t-x_j\,|\,j\ne i]_l$). Indeed, assuming
that the $g_1(x_1)\ne 0,\dots,g_n(x_n)\ne 0$, we get that the
matrix $C$ in formula (F)(just before \ref{QD of VDM}) is
invertible. On the other hand, for any $i=1,\dots,n \;g_i(x_i)\ne
0$ implies that
$deg(g_i)+1=deg([g_i(t),t-x_i]_l)=deg([t-x_1,\dots,t-x_n]_l)<n$,
by our assumption. Hence $deg(g_i)\le n-2$, this implies that the
last column of the matrix $C$ is zero.   This contradicts the
invertibility of $C$ and proves the claim.  Assume thus that there
exists $1\le i\le n$ such that $g_i(x_i)=0$.  Comparing the $ith$
row on both sides of the equation (F) we get
$(c_{i,0}\dots,c_{i,n-1})V_n^{S,D}(x_1,\dots,x_n)=(0,\dots,0)$.
Since $g_i\ne 0$ and $V_n^{S,D}(x_1,\dots,x_n)$ is invertible this
gives a contradiction.

\noindent  $i) \Leftrightarrow i')$.  It is clear that $i')$
implies $i)$.  For the converse implication, let us put
$p_n(t):=[t-x_i\,|\, 1\le i\le n ]_l$ and suppose that
$deg(p_n)=n$.  Define $A:=\{i_1,\dots,i_r\}\subseteq
\{1,\dots,n\}$, $A':=\{1,\dots,n\}\setminus A$, $p_A(t):=[t-x_i \,
|\,i\in A]_l$ and $p_{A'}(t):=[t-x_j \, | \, j\in A']_l$.  We then
have $p_n(t)=[t-x_i \, | \, i\in A \cup
A']_l=[p_A(t),p_{A'}(t)]_l$ and since $deg(p_n(t))=n$, we easily
get that $deg(p_A(t))=|A|=r$, as desired.

\noindent $i') \Leftrightarrow ii')$.  This is clear from $i)
\Leftrightarrow ii)$.

\noindent Let us now suppose that $x_1=a^{u_1},\dots,x_n=a^{u_n}$.

\noindent $i) \Fl iii)$.  Assume $i)$ holds but $u_1,\dots,u_n$
are right $C^{S,D}(a)$ linearly dependent.  Without loss of
generality we may assume that $u_n=\sum_{i=1}^{n-1}u_ic_i$ where
$c_1,\dots,c_{n-1} \in C^{S,D}(a)$.   Let us put
$p_{n-1}(t)=[t-x_i \,|\,i=1,\dots,n-1]_l$.  Using Lemma
\ref{properties of T_a} we obtain
$p_{n-1}(a^{u_n})=p_{n-1}(T_a)(u_n)u_n^{-1}=
p_{n-1}(T_a)(\sum_{i=1}^{n-1}u_ic_i)u_n^{-1}
=\sum_{i=1}^{n-1}p_{n-1}(T_a)(u_i)c_iu_n^{-1}=
\sum_{i=1}^{n-1}p_{n-1}(a^{u_i})u_i^{-1}c_iu_n^{-1}=0$ where the
last equality comes from the definition of $p_{n-1}(t)$.  This
shows that $p_{n-1}(x_{n})=0$ and contradicts $i)$.

\noindent $iii) \Fl i)$.  We prove, by induction on $n$, that
$Ker(p_n(T_a) )=\sum_{i=1}^nu_iC^{S,D}(a)$ and that
$deg(p_n(t))=n$.   If $n=1$, $p_1(t)=t-a^{u_1}$ is of degree $1$
and $p_1(T_a)(v)=0$ implies that either $v=0$ or $T_a(v)=a^{u_1}v$
which leads to $a^{u_1}=a^v$ and hence $v\in u_1C^{S,D}(a)$.

\noindent Suppose $n>1$ and assume, by induction, that
$Ker(p_{n-1}(T_a) )=\sum_{i=1}^{n-1}u_iC^{S,D}(a)$ and
$deg(p_{n-1}(t) )=n-1$.  If $p_n(t)=p_{n-1}(t)$ we get $u_n\in
Ker(p_n(T_a) )=Ker(p_{n-1}(T_a) )=\sum_{i=1}^{n-1}u_iC^{S,D}$, a
contradiction.  Hence we  must have $p_n(t)\ne p_{n-1}(t)$.  Using
lemma \ref{product formula} 2. and \ref{properties of T_a} one can
write
$p_n(t)=[t-x_n,p_{n-1}(t)]_l=(t-x_n^{p_{n-1}(x_n)})p_{n-1}(t)=(t-a^{w_n})p_{n-1}(t)$,
where $w_n=p_{n-1}(a^{u_n})u_n=(p_{n-1}(T_a) )(u_n)$.  In
particular, $deg p_n(t)=n$.  The induction hypothesis shows that
$Ker(p_{n-1}(T_a))=\sum_{i=1}^{n_1}u_iC^{S,D}(a)\subseteq
\sum_{i=1}^{n}u_iC^{S,D}(a)\subseteq Ker(p_n(Ta) )$.  Moreover, if
$v\in Ker(p_n(T_a) )\setminus Ker(p_{n-1}(T_a)
)=\sum_{i=1}^{n-1}u_iC^{S,D}(a)$, we have $T_a(p_{n-1}(T_a)(v)
)=a^{w_n}p_n(T_a)(v)$.  This leads to
$p_{n-1}(T_a)(v)=w_nc=p_{n-1}(T_a)(u_nc)$.  Hence $v-u_nc \in
Ker(p_{n-1}(T_a) )=\sum_{i=1}^{n-1}u_iC^{S,D}(a)$.   This yields
the conclusion.


\noindent $ii) \Leftrightarrow iv)$.  This is extracted from
Proposition \ref{QD of WRSK}.

\noindent The other implications are clear.
\end{proof}

Let us recall that for $f\in R=K[t;S,D]$ and $a\in K$, $E(f,a)$
stands for the set $\{x\in K\, |\, f(a^x)=0\}\cup\{0\}$.  Lemma
\ref{properties of T_a} shows that $E(f,a)=Ker(f(T_a))$.

\begin{corollary}
\label{dim E(f) and deg(f)}
\begin{enumerate}
\item For $f\in R$ and $a\in K$, we have $dim_{C(a)}E(f,a)\le
deg(f)$, where $C(a)$ stands for $C^{S,D}(a)$.
\item If $u_1,\dots,u_n\in K$ are right $C^{S,D}(a)$-linearly
independent and $p_n(t)=[t-a^{u_i}\,|\, 1\le i \le n]_l$ then
$Ker(p_n(T_a) )=\sum_{i=1}^nu_iC^{S,D}(a)$.
\end{enumerate}
\end{corollary}

\begin{proof}
1. Assume, at the contrary, that $dim_{C(a)}E(f,a) > deg(f)=:n$.
This means that there exist $u_1,\dots,u_{n+1}\in E(f,a)$ which
are linearly independent over $C(a)$.  Let us put
$g:=[t-a^{u_1},\dots,t-a^{u_{n+1} }]_l$; $g$ is a W-polynomial of
degree $n+1$. Since $f(a^{u_i})=0$ for $i=1,\dots,n+1$, we have
$Rf \subseteq \cap_{i=1}^{n+1}R(t-a^{u_i})=Rg$.  Hence $deg(f)\ge
n+1$.  This contradiction yields the lemma.

2.  This has been shown in the proof of the implication $iii) \Fl
i)$ in the above Theorem \ref{equivalence for P-independence}.
\end{proof}

\begin{remarks}
{\rm
\begin{enumerate}
\label{remark on P-independence}
\item[a)] The notion of $P$-independence comes from an abstract
dependence relation.  Let us recall the definition :
\begin{definition}
Let $S$ be a set.  A dependence relation on $S$ is a rule which
associates with each finite subset $X$ of $S$ certain elements of
$S$, said to be dependent on $X$.  The following conditions must
be satisfied :
\begin{enumerate}
\item If $X=\{x_1,\dots,x_n\}$, then each $x_i$ is dependent on
$X$.
\item If $z$ is dependent on $\{y_1,\dots,y_n\}$ and each $y_i$ is
dependent on $\{x_1,\dots,x_n\}$, then $z$ is dependent on
$\{x_1,\dots,x_n\}$.
\item If $y$ is dependent on $\{x_1,\dots,x_n\}$, but not on
$\{x_2,\dots,x_n\}$, then $x_1$ is dependent on
$\{y,x_2,\dots,x_n\}$.
\end{enumerate}
\end{definition}
In our situation an element $y\in K$ is $P$-dependent on a finite
subset $X$ of the division ring $K$ if $[t-x\,|\, x\in X]_l(y)=0$.
We leave to the reader the easy proof that this defines indeed a
 dependence relation on $K$.  Actually these notions can be put in
the more general frame of $2$-firs (Cf. \cite{LO}).

\item[b)] In \cite{GRW} section $4.1$ the Miura decomposition of
differential operators $L(D)=D^n+a_{n-1}D^{n-1}+\dots+a_1D+a_0$
with coefficients $a_i$ in a division ring $K$ was given.  In this
paper the authors assumed that :
\begin{enumerate}
\item There exist $n$ independent solutions say $u_1,\dots,u_n$ of the equation
$l(D)=0$ (over the subdivision ring $Ker D$) ).
\item They also assumed that for any subset $\{\i_1,\dots,i_r\}\subseteq \{1,\dots,n\}$
the Wronskian matrix $W(u_{i_1},\dots,u_{i_r} )$ is invertible.
\end{enumerate}
Lemma \ref{properties of T_a} shows that for $0\ne x\in K$ and
$L(t)\in K[t;Id.D]$,\; $L(D)(x)=0$ if and only is $L(0^x)=0$.
Notice also that $Ker D=C^{Id.,D}(0)$.  The equivalence between
$iii)$ and $iv')$ shows that this second hypothesis is uncessary.

\end{enumerate}
}
\end{remarks}

As mentioned in remark \ref{LU in general U algorithmically} b),
it is worth to notice that the matrix $U$ in \ref{LU decomposition
of V and W} can be algorithmically computed.  This is the aim of
the next proposition.   Moreover it offers a way of testing if a
given set $\{x_1,\dots,x_n\}\subseteq K$ is $P$-independent.  The
diagonal elements of the matrix $U$ also gives an algorithm for
computing a linear factorization of a $W$-polynomial
$[t-x_j\,|\,j=1,\dots,n]$.

\vspace{3mm}

\begin{prop}
\label{algorithm for computing U} Let $x_1,\dots,x_n$ be any
finite set of elements of $K$ and define inductively the elements
$u_{ij}\in K$ for $i=0,\dots,n-1$ and $j=1,\dots,n$ as follows:
$u_{0j}=1$ for $j=1,\dots,n$ and assuming that
$u_{i,1},\dots,u_{i,n}$ have been defined we put
$$
u_{i+1,j}=\left \lbrace
                 \begin{array}{l}
                \quad \quad \quad 0 \quad \quad \quad \quad \; if \quad u_{i,j}u_{i,i+1} =0\\
                 (x_{j}^{u_{ij} }-x_{i+1}^{u_{i,i+1} })u_{ij} \quad if \quad
                 u_{ij}\ne 0\ne u_{i,i+1}.
                 \end{array}
                 \right. \quad.
$$
The following are equivalent :
\begin{enumerate}
\item[a)] The set $\{x_1,\dots,x_n\}$ is $P$-independent.
\item[b)] $u_{ij}\ne 0 $ for all $(i,j)\in \mathbb N$ such that $0\le i < j \le n$.
\item[c)] $u_{n-1,n}\ne 0$.
\end{enumerate}
In this case, the matrix $U=(u_{ij})$ is the one obtained in
Proposition \ref{LU decomposition of V and W} and one has :
$[t-x_j\,|\, j=1,\dots,n]_l=(t-x_n^{u_{n-1,n} })\cdots
(t-x_2^{u_{12} })(t-x_1^{u_{01} })$.

\end{prop}
\begin{proof}
Let us define, for $1\le i\le n,\;
p_i(t):=[t-x_j\,|\,j=1,\dots,i]_l$ and $p_o(t)=1$.

\noindent $a) \Longrightarrow b)$. Suppose that
$\{x_1,\dots,x_n\}\subseteq K$ is a $P$-independent set, then
$p_{i+1}(t)=[t-x_{i+1},p_i(t)]_l=(t-x_{i+1}^{p_i(x_{i+1})})p_i(t)$.
Using this formula it is easy to prove, by induction on $i$ $0\le
i\le n-1$, that $u_{ij}=p_i(x_j)$.   The independence of the set
$\{x_1,x_2\dots,x_n\}$ yields that $u_{ij}\ne 0$ for $j>i$.

\noindent b) $\Fl $ c) This is clear.

\noindent c) $\Fl$ a) Suppose that $u_{n-1,n}\ne 0$.  The
definition of this element shows that $u_{n-2,n}\ne 0$ and
$u_{n-2,n-1}\ne 0$.  Continuing this process "backwards" we get
that $u_{ij}\ne 0$ for all $(i,j)$ such that $0\le i<j\le n$.  Now
let us show, by induction on $i$, that for $0\le i<j\le n$ we have
$u_{ij}=p_i(x_j)$.  For $i=0$, we have $u_{0j}=1=p_0(x_j)$ (recall
that $P_0(t)=1$).  Assume we have proved $u_{ij}=p_{i}(x_j)$ for
all $j>i$ and let us consider $u_{i+1,j}$ for $j>i+1$.  We have :

\noindent $u_{i+1,j}=(x_j^{u_{ij} }-x_{i+1}^{u_{i,i+1} })u_{ij}=
                (x_j^{p_i(x_j)}-x_{i+1}^{p_i(x_{i+1}) })p_i(x_j)
=((t-x_{i+1}^{p_i(x_{i+1})})p_i(t))(x_j)\\
= [t-x_{i+1},p_i(t)]_l(x_j)=p_{i+1}(x_j)$.  This ends the
induction and shows that $p_{i+1}(x_j)=u_{i+1,j}\ne 0$, for
$j>i+1$.   We conclude that $p_i(x_j)\ne 0$ for $0\le i<j\le n$.
In particular, $\{x_1,\dots,x_n\}$ is a $P$-independent set.

The other statements are now clear.
 \end{proof}

\vspace{10mm}

\section{Linear factorizations of W-polynomials.}

Let us introduce some notations  and a definition:

\noindent For $f\in R=K[t;S,D]$ we denote by $V(f)$ the set of
right roots of $f$ i.e. $V(f):=\{a\in K \, | \, f(a)=0\}=\{a \in K
\,|\, f(t)\in R(t-a)\}$.  In case $S=Id.$ and $D=0$ it was proved
by Gordon and Motzkin (Cf. \cite{GM}) that $V(f)$ intersects a
finite number of conjugacy classes.  This is also true in an
$(S,D)$-setting (Cf. \cite{LLO}).  We can thus write
$V(f)=\cup_{i=1}^rV_i$ where $V_i=V(f)\cap \De^{S,D}(a_i)$. For
$i=1,\dots,r$ the set $E(f,a_i):=\{x \in K \setminus {0}\,|\,
a_i^x \in V_i \}\cup \{0\}$ is in fact a right vector space over
the division ring $C_i=C^{S,D}(a_i):=\{x \in K \setminus\{0\}
\,|\, a_i^x=a_i \}\cup \{0\}$.  It was proved in $\cite{LLO}$ that
:
$$
\sum_{i=1}^r dim_{C_i}E(f,a_i) \le deg(f)
$$
Moreover, in this formula the equality holds if and only if $f$ is
a W-polynomial.

It is easy to remark that for any $f\in R$, there exists a
W-polynomial $g\in R$ such that $V(f)=V(g)$.  For the rest of this
section $f$ will stand for a W-polynomial of degree $n$.

\begin{definition}
Let $f\in R=K[t;S,D]$ be a W-polynomial of degree $n$.  A
$P$-basis for $V(f):=\{x\in K\,|\, f(x)=0\}$ is a set
$\{x_1,\dots,x_n\}\subseteq V(f)$ such that $[t-x_i \,|\,
i=1,\dots,n]_l=f$.
\end{definition}

We are interested in describing all the different linear
factorizations of $f$.  Let us first consider the case when all
the roots of $f$ are in a single conjugacy class :
$\Delta^{S,D}(a)$.  If $x_1=a^{u_1},\dots,x_n=a^{u_n}$ is a
$P$-basis for $V(f)$, Theorem \ref{equivalence for P-independence}
shows that all the $P$-bases of $f$ are of the form
$a^{v_1},\dots,a^{v_n}$ where $(v_1,\dots,v_n)=(u_1,\dots,u_n)A$
for some matrix $A\in GL_n(C^{S,D}(a))$.  To every ordered P-basis
$(x_1,\dots,x_n)$ we can associate, as in section $2$, a
factorization of $f$ given by
$f(t)=(t-x_n^{z_n})\cdots(t-x_i^{z_i})\cdots (t-z_1)$ where
$z_i=p_{i-1}(x_i)$ and $p_i=[t-x_1,\dots,t-x_i]$ for
$i=1,\dots,n-1$.  But different ordered P-bases can lead to the
same factorization.  Let us give an easy example :

\begin{example}
\label{different P bases but same fact} {\rm Let $K$ be a division
ring $S=Id.$ and $D=0$.  Suppose that $a\in K$ is such that $dim
K_{C(a)}\ge 3$ where $C(a)$ denotes the usual centralizer of $a$.
Let $u_1,u_2,u_3$ be three elements in $K$ which are right
linearly independent over $C(a)$ and consider
$f=[t-a^{u_i}\,|\,i=1,2,3]_l$.  The elements
$v_1=u_1,v_2=u_1+u_2,v_3=u_3$ form another right basis of $K$ over
$C(a)$.  The two ordered bases $(a^{u_1},a^{u_2},a^{u_3})$ and
$(a^{v_1},a^{v_2},a^{v_3})$ are different but give rise to the
same factorization.  Indeed, putting $x_i:=a^{u_i}$ and
$y_i:=a^{v_i}$ for $1\le i \le 3$, we have $p_1(t)=t-x_1=t-y_1$,
$x_2^{p_1(x_2)}=a^{p_1(x_2)u_2}$ and
$y_2^{p_1(y_2)}=a^{p_1(y_2)v_2}$; but $p_1(y_2)v_2=p_1(a^{v_2})v_2
=p_1(T_a)(v_2)=p_1(T_a)(u_1+u_2)=p_1(T_a)(u_2)=p_1(a^{u_2})u_2=p_1(x_2)u_2$
and we conclude that $x_2^{p_1(x_2)}=y_2^{p_1(y_2)}$.  It is then
easy to check that the two factorizations given by the different
ordered $P$-bases are the same.  Details are left to the reader. }
\end{example}

The next lemma will be very useful.

\begin{lemma}
\label{a^g(a) and b^g(b)} \hfill
\begin{enumerate}
\item[a)] $a^{p(T_a)(v)}=a^{p(T_a)(u)}$ if and only if $Ker p(T_a)+vC=Ker
p(T_a)+uC$.
\item[b)] Let $g\in R$ be any
polynomial and $a,b\in K$ be distinct elements of $K$ such that
$g(a)\ne 0$ and $g(b)\ne 0$. Then $a^{g(a)}=b^{g(b)}$ if and only
if there exist $u\in K\setminus C^{S,D}(a)$ and $c\in C^{S,D}(a)$
such that $b=a^u$ and $g(a^{u-c})=0$.  In particular, in these
circumstances we have $V(g)\cap \Delta(a)\ne \emptyset$.
\end{enumerate}
\end{lemma}
\begin{proof}
\noindent a) We have $a^{p(T_a)(v)}=a^{p(T_a)(u)}$ if and only if
there exists $c\in C^{S,D}(a)$ such that $p(T_a)(v)=p(T_a)(u)c$.
Since $p(T_a)$ is a right $C^{S,D}(a)$-linear map this is
equivalent to $p(T_a)(v-uc)=0$.  This easily gives the result.

\noindent b) Since $a^{g(a)}=b^{g(b)}$, we can write $b=a^u$ for
some $u\in K\setminus C^{S,D}(a)$.  We then have
$a^{g(T_a)(1)}=a^{g(a)}=b^{g(b)}=a^{g(T_a)(u)}$, hence there
exists $c\in C^{S,D}(a)$ such that
$g(T_a)(u)=g(T_a)(1)c=g(T_a)(c)$ and $u-c \in Ker(g(T_a))$ i.e.
$g(a^{u-c})=0$.  The converse implication is left to the reader.
\end{proof}

\begin{theorem}
\label{bijection flags in E and factorizations when all rooots in
1 CC } Let $f\in R=K[t;S,D]$ be a W-polynomial of degree $n$ such
that $V(f)\subseteq \De^{S,D}(a)$.  The linear factorizations of
$f$ are in bijection with the complete flags of the right
$C^{S,D}(a)$ vector space $E(f,a)$.
%

\end{theorem}
\begin{proof}
%
Let $f(t)=(t-a_n)\cdots(t-a_1)$ be a linear factorization of the
$W$-polynomial $f$.  The polynomial
$p_i(t):=(t-a_i)(t-a_{i-1})\dots(t-a_1)$ is also a $W$-polynomial
and the fact that $V(f)\subseteq \Delta^{S,D}(a)$ implies that we
can write $Kerp_i(T_a)=\sum u_jC^{S,D}(a)$ for some right
$C^{S,D}(a)$-independent elements $u_1,\dots,u_i$ in
$E(p_i,a)\subseteq E(f,a)$ (Cf. Corollary \ref{dim E(f) and
deg(f)} and \ref{characterization of W-polynomials}).

\noindent Let $\psi$ be the map from the set of factorizations of
$f$ to the set of flags in $E(f,a)$ defined by associating the
flag $Ker(p_1(T_a))\subset Ker(p_2(T_a))\subset \dots \subset
Ker(p_n(T_a))=E(f,a)$ to the factorization
$f=(t-a_n)\cdots(t-a_1)$.

\noindent Let us show that $\psi$ is injective : if
$f(t)=(t-a_n)\cdots(t-a_1)=(t-a'_n)\cdots(t-a'_1)$ are two
different factorizations of $f$, we define $l:=min \{i \,|\,
a_i\ne a'_i \}$.  We thus have
$p_i(t)=(t-a_i)\cdots(t-a_1)=(t-a'_i)\cdots(t-a'_1)$ for any
$i<l$, $p_l(t)=(t-a_l)p_{l-1}(t)$ and $p'_l(t)=(t-a'_l)p_{l-1}(t)$
(if $l=1, \; p_{l-1}(t)=p_o(t)=1$).  Notice that $p_l(t)$ and
$p'_l(t)$ are W-polynomials.  Hence one can write
$p_l(t)=[t-a^u,p_{l-1}(t)]$ and $p'_l(t)=[t-a^v,p_{l-1}(t)]$ for
some $u,v\in E(f,a)$.  This gives $a_l=a^{p_{l-1}(T_a)(u)}$ and
$a'_l=a^{p_{l-1}(T_a)(v)}$.  Lemma \ref{a^g(a) and b^g(b)} a)
shows that the flags associated with these two factorizations will
be different.

\noindent Let us show that $\psi$ is onto.  Consider a complete
flag $u_1C\subset u_1C+u_2C \subset \dots \subset
\sum_{i=1}^nu_iC=E(f,a)$.  We build successively the following
right factors of $f$ :
$p_o(t)=1,p_1(t)=t-a_1,\dots,p_i(t)=[t-a^{u_i},p_{i-1}(t)]_l$.
This will give a factorization of $f=p_n$ :
$f(t)=(t-a^{p_{n-1}(T_a)(u_n)}(t-a^{p_{n-2}(T_a)(u_{n-1})})\dots
(t-a^{u_1})$.  It is then easy to check that $\psi$ maps this
factorization on the complete flag we started with.

\end{proof}

\begin{example}
\label{factorizations in dimension 2} {\rm Let us describe all the
factorizations of $f=[t-a^x,t-a]_l$.  These factorizations are in
bijection with the complete flags in the two dimensional vector
space $E(f,a)=C+xC$ where $C:=C^{S,D}(a)$.  These flags are of the
form $0\ne yC \subset E(f,a)$.  Apart from the flag $0 \subset xC
\subset E(f,a)$, they are given by $0 \subset (1+x\beta)C \subset
E(f,a)$, where $\beta \in C^{S,D}(a)$.  Hence we get the following
factorizations $f=(t-a^{a-a^x})(t-a^x)$ and
$(t-a^{a-\gamma})(t-a^{1+x\beta})$, where $\gamma=a-a^{1+x\beta}$.

}
\end{example}

Let us now describe all linear factorizations of a general
$W$-polynomial $f$ (i.e. without assuming that all its roots are
in a single conjugacy class).  Before stating our last theorem in
this section let us fix some notations.  For a $W$-polynomial $f$
we decompose $V(f)$ with respect to conjugacy classes :
$V(f)=\bigcup_{i=1}^rV_i$ where $V_i=V(f)\cap \De(a_i)$.  Consider
the semisimple ring $C:=\prod_{i=1}^rC_i$ where, for
$i=1,\dots,r$, $C_i:=C^{S,D}(a_i)$.   $K^r$ has a natural
structure of right $C$-module. Its submodules are all of the form
$U_1\times \cdots \times U_r$, where, for $i=1,\dots,r$,
$U_i\subseteq K$ is a right $C_i$ vector space.  In particular,
$E(f):=\prod_{i=1}^rE(f,a_i)$ is a right $C$-module.  Notice that
the dimensions $dim_{C_i}E(f,a_i)<\infty$.

\begin{definitions}  Let $f$ be a W-polynomial in $R=K[t;S,D]$.
\begin{enumerate}
\item[a)] For a $C$-submodule $U= U_1\times U_2\times\dots\times U_r$ of $E(f)$ we
define its dimension to be
$dim(U):=(dim_{C_1}U_1,\dots,dim_{C_r}U_r)$ and its weight to be
$wt(U):=\sum_{i=1}^rdim_{C_i}U_i$.
\item[b)] A complete flag in $E(f)$ is a sequence of $C$-submodules
$0\subsetneqq M_1 \subsetneqq M_2 \subsetneqq \dots \subsetneqq
M_n$ such that $n=deg(f)$.
\item[c)] For $a=(a_1,\dots,a_r)\in K^r$ we define $T_a : K^r \fl K^r:(x_1,\dots,x_r)\mapsto
(T_{a_1}(x_1),\dots,T_{a_r}(x_r))$.
\end{enumerate}
\end{definitions}

Let us notice that, by the remark stated at the end of the first
paragraph introducing this section, if $f$ is a $W$-polynomial the
weight of the right $C$-module $E(f)$ is $wt(E(f))=deg(f)$. Notice
also that a sequence of right $C$-modules $0\subsetneqq M_1
\subsetneqq M_2 \subsetneqq \dots \subsetneqq M_n$ is a complete
flag if and only if $wt(M_i)=i$.

\noindent With these notations we can state the next result which
generalizes Theorem \ref{bijection flags in E and factorizations
when all rooots in 1 CC }.

\begin{theorem}
\label{bijection factorizations and flags general case} Let $f$ be
a W-polynomial. Then :

\noindent There is a bijection between the set of factorizations
of $f$ and the set of complete flags in $E(f)=\prod_{i=1}^rE_i$.
\end{theorem}
\begin{proof}
As in the proof of \ref{bijection flags in E and factorizations
when all rooots in 1 CC } we associate to the factorization
$f(t)=(t-b_1)\dots (t-b_n)$ the flag of right $C$-module given by
$0\subsetneqq Ker(p_1(T_a))\subsetneqq \dots \subsetneqq
Ker(p_n(T_a))$ where $p_i(t)=(t-b_{i})(t-b_{i-1})\cdots(t-b_1)$.
Let us show, by induction on $i$, that $wt(Ker(p_i(T_a) ) )=i$. If
$i=1$ we have $(0,\dots,0)\ne (x_1,\dots,x_n)\in Ker(p_1(T_a))$ if
and only if $(0,\dots,0)=(T_a-b_1)(x_1,\dots,x_n)=(S(x_1)a_1
+D(x_1)-b_1x_1, \dots,S(x_r)a_r +D(x_r)-b_1x_r)$.  This shows that
for $1\le i \le r$, we have $S(x_i)a_i + D(x_i)=b_1x_i$. We can
assume that $b_1\in \Delta(a_1)$, say $b_1=a_1^{v_1}$.  Now assume
that for some $j>1$, $x_j\ne 0$.  we then have
$b_1=S(x_j)a_jx_j^{-1} + D(x_j)x_j^{-1}$.  In particular $b_1\in
\Delta^{S,D}(a_j)$.  This contradicts the fact that $b_1 \in
\Delta (a_1)$.  We conclude that $Ker(p_1(T_a))=v_1C^{S,D}(a_1)$,
as desired.  Assume that $wt(Ker(p_i(T_a) ) )=i$ and let us show
that $wt(Ker(p_{i+1}(T_a) ) )=i+1$.  We have
$p_{i+1}=(t-b_{i+1})p_i(t)$.  If $v\in Ker(p_{i+1}(T_a))\setminus
Ker(p_i(T_a))$ and $b_{i+1}\in \Delta (a_j)$, say
$b_{i+1}=a_j^{v_j}$, for some $j\in\{1,\dots,r\}$ and some $v_j\in
K\setminus\{0\}$.  As in case $i=1$ we easily check that
$v=(0,\dots,x_j,0\dots,0)$ where $x_j=v_jC^{s,D}(a_j)$.  This
shows that $wt(Ker(p_{i+1}(T_a) ) )=wt (Kerp(p_i(T_a)))+1$.  The
induction hypothesis implies that $wt(Ker(p_{i+1}(T_a) ) )=i+1$,
as desired.   The rest of the proof is completely similar to the
one given in \ref{bijection flags in E and factorizations when all
rooots in 1 CC } abd is left to the reader.

\end{proof}
\begin{remark}
\begin{enumerate}
\item[a)] Let us remark that $\De(a)$ has a structure of right
$C^{S,D}(a)$ projective space which is in fact given by the right
$C^{S,D}(a)$ vector space structure of $K$ itself (the map $\phi :
\mathbb P(K_{C(a)})\longrightarrow \De(a):x \mapsto a^x$ is a
bijection). If $V(f)\subseteq  \De(a)$ then $V(f)$ is a projective
subspace of $\mathbb P(K)$, the associated vector space being
$E(f,a)$.
\item[b)]Gelfand, Retakh and Wilson introduced and studied in
details an algebra $Q_n$ associated with factorizations of certain
polynomials $f(t)$ in the universal field of quotients $K$ (aka.
free field) of the free algebra $k<x_1,x_2,\dots,x_n>$ over a
commutative field $k$.  The aim is to replace te free field by
some "smaller" algebra in which all the factorizations of $f(t)$
already take place.  Using the least left common multiple, the
introduction and the description of this algebra are very natural.
In this language, the polynomial of which we study the
factorization is $f(t)=[t-x_i\,|\, i\{1,\dots,n\}]_l\in K[t]$.
Let us first fix some notations.  For $A\subseteq \{1,\dots,n\}$
we denote $p_A:=[t-x_i\,|\, i\in A]_l\in K[t]$.
  Let $i,j\in\{1,\dots,n\}\setminus A$ we have
$p_{A\cup \{i,j\} }(t)=[t-x_j,[t-x_i,p_A(t)]_l ]_l=
[t-x_i,[t-x_j,p_A(t)]_l]_l$. Putting $x_{A,i}:=x_i^{p_A(x_i)}$ and
using \ref{product formula} $2.$ we obtain :

\noindent $(t-x_{A\cup \{i\},j} )(t-x_{A,i})p_A(t)=
 (t-x_{A\cup \{j\},i})(t-x_{A,i})p_A(t)$. This leads to the quadratic relations :
$$
x_{A\cup\{j\},i} +x_{A,j}=x_{A\cup\{i\},j} +x_{A,i} \quad{\rm
and}\quad x_{A\cup\{j\},i}.x_{A,j}=x_{A\cup\{i\},j}.x_{A,i} .
$$

The algebra is then describe by $Q_n:=k<z_{A,i} \,|\, A\subseteq
\{1,\dots,n\}>/I $, where $I$ is the ideal generated by the
analogue of the quadratic relations obtained above.  All the
factorizations we have obtained are also obtainable in $Q_n[t]$.
This algebra $Q_n$ have been studied in depth by Gelfand, Retakh
and Wilson (Cf.\cite{GRW}).  Notice also that it is possible to
formally introduce more generally such an algebra in an
$(S,D)$-setting.  This might be an interesting way of looking to
factorizations of differential operators.
\end{enumerate}
\end{remark}

\section{Existence of LLCM}

The aim of this short section is to have a nice criterion for a
ring $A$ to be such that for any finite subset
$\{a_1,\dots,a_n\}\subseteq A$, there exists a monic polynomial
$p(t)\in A[t],\; deg(p(t))=n$ such that, for any $i\in
\{1,\dots,n\}$, $t-a_i$ divides $p(t)$ on the right.  With our
standard notation we can then write $p(t)=[t-a_i \, | \,
i=1,\dots,n]$.  Recall that a ring is called left duo if its left
ideals are in fact two-sided ideals.


\begin{lemma}
\label{llcm for 2 gives llcn for finite} Let $A,S,D$ be a ring an
endomorphism of $A$ and an $S$-derivation of $A$, respectively.
Assume that for any $a,b\in A$ there exist $c,d \in A$ such that
$(t-c)(t-a)=(t-d)(t-b)\in R=A[t;S,D]$ then for any finite subset
$\{a_1,\dots,a_n\}\subseteq A$ there exists a monic polynomial of
degree $n$ in $\cap_{i=1}^nR(t-a_i)$.
\end{lemma}
\begin{proof}
We first construct, for any $1\le l \le n $, elements $a_{i_1\dots
i_l}$  with $1\le i_j \le n$ for all $j=1\dots l$ and $i_j\ne i_k$
if $j\ne k$.

\noindent If $l=1$, $a_{i_1}$ is given and belongs to $\{a_1,\dots
a_n\}$.

\noindent For $l\ge 2$ we assume that the elements $a_{i_1\dots
i_{l-1}}\in A$ have been constructed for any set of distinct
indexes $\{ i_1,\dots,i_{l-1}\}$ with $1\le i_j \le n$ for any
$1\le j \le l-1$.  By hypothesis, there exist elements $a_{i_1
\dots i_{l-1}i_l}$ and
 $a_{i_1 \dots
i_{l-2}i_li_{l-1} }$ in $A$ such that :

$$(t-a_{i_1\dots i_l})(t-a_{i_1\dots i_{l-1} })= (t-a_{i_1 \dots
i_{l-2}i_li_{l-1} })(t-a_{i_1\dots i_{l-2}i_{l} }).$$

This defines $a_{i_1 \dots i_{l-1}i_l}$ for any subset
$\{i_1,\dots,i_l\}$ of $\{1,\dots,n\}$.  Let us put
$f_1(t):=t-a_1$ and $f_l(t):=(t-a_{1\dots l})f_{l-1}(t)$, for $1<
l < n$.  First let us show, by induction on $l$, that for any $l <
s \le n$ we have $(t-a_{1\dots ls})f_l \in R(t-a_s)$. For $l=1$,
$f_1=t-a_1$ and we have $(t-a_{1s})(t-a_1)=(t-a_{s1})(t-a_s)\in
R(t-a_s)$.  For $l>1$, we have $(t-a_{1\dots ls})f_l= (t-
a_{1\dots ls})(t-a_{1\dots l})f_{l-1}=(t-a_{1 \dots
(l-1)sl})(t-a_{1\dots (l-1)s})f_{l-1}\in R(t-a_s)$ by the
induction hypothesis.

\noindent We can now prove easily, by induction on $l$, that
$f_l\in \cap_{i=1}^lR(t-a_i)$.  This is left to the reader.
 $f_n$ is then the required monic polynomial of degree $n$ in
$\cap_{i=1}^nR(t-a_i)$
\end{proof}

\begin{theorem}
\label{existence of llcm} Let $A,S,D$ be a ring an endomorphism of
$A$ and a $S$-derivation of $A$, respectively.  For $R=A[t;S,D]$
the following are equivalent :
\begin{enumerate}
\item For any $\{a_1,\dots,a_n\}\subseteq A$ there exists a monic
polynomial of degree $n$ in $\cap_{i=1}^nR(t- a_i) \subseteq
R[t]$.
\item For any $\{a,b\}\subseteq A$ there exists a monic polynomial of degree
two in $R(t-a)\cap R(t-b)$.
\item For any $r,b\in A$ there exists $c\in A$ such that $cr=S(r)b
+D(r)$.
\end{enumerate}
\noindent If $S=Id.$ and $D=0$ the above statements are also
equivalent to $A$ being left duo.

\end{theorem}
\begin{proof}

1) $\Leftrightarrow $ 2) This is clear thanks to the above lemma
\ref{llcm for 2 gives llcn for finite}.

2) $\Fl$ 3) Let $r,b\in A$.   By hypothesis we know that there
exists $c,d\in A$ such that $(t-c)(t- (r+b))=(t-d)(t-b)\in R=
A[t;S,D]$. Equating coefficients of the same degree gives
$c(r+b)-D(r+b)=db-D(b)$ and $c+S(r+b)=d+S(b)$. Replacing $d$ in
the first equality we get $cr-D(r)=S(r)b$, as desired.

3) $\Fl$ 2) Assume that for any $r,b\in A$ there exists $c\in A$
such that $cr-D(r)=S(r)b$ and let $a,b\in A$ be given.  There
exists $c\in A$ such that $c(a-b) =S(a-b)b+D(a-b)$. Putting
$d:=c+S(a-b)$ we get that $(t-d)(t-b)=(t-c)(t-a)\in A[t;S,D]$, as
required.

\noindent Now, if $S=id.$ and $D=0$, 3. above is equivalent to the
fact that, for any $r,b\in A$, $Ab\in Ar$.  This means that $Ar$
is a two-sided ideal.  This means that every principal left ideal
is in fact two-sided i.e. $A$ is left duo.
\end{proof}

\begin{remark}
The last statement of the above theorem shows that in an $S,D$
setting the equality in \ref{existence of llcm} 3. could be
considered as a definition for an $(S,D)$ left duo ring.
\end{remark}

\begin{example}
{\rm Let $k$ be a field.  Obviously $M_2(k)$ is not left duo.
Indeed it is easy to remark that the matrices
$$
              a=\begin{pmatrix}
                 1   &  1  \\
                 0   &  1
                 \end{pmatrix}\quad
{\mbox {\rm and } }\quad
              b=\begin{pmatrix}
                0    &  0  \\
                0    &  1
              \end{pmatrix}
$$
are such that $t-a$ and $t-b$ have no left monic common multiple
of degree two.   Of course this also means that there are no
matrices $c,d\in M_2(k)$ such that $c+a=d+b$ and $ca=db$.  This
also gives that there is no solution to the matrix equation $t^2 -
(a+b)t + ab =0$.
}
\end{example}

We exhibit now an example showing that the left and right
existence of a least monic common multiple are different notions.

\begin{example}
{\rm Consider $A$ the subring of upper triangular matrices over
$\mathbb Q (x)$ define as follows :
$$
A=\{ \begin{pmatrix}
                 f(x^2) & g(x) \\
                   0    & f(x)
     \end{pmatrix} \,|\, f(x),g(x)\in \mathbb Q (x)\}
$$
Let $R$ denote the usual polynomial rings $R=A[t]$.  It is easy to
check that there is no monic polynomial of degree $2$ in
$R(t-a)\cap R(t-b)$ where
$$
a=\begin{pmatrix}
                 x^2 & x^2 \\
                   0 &  x
     \end{pmatrix} \quad {\mbox {\rm and } }\quad b=\begin{pmatrix}
                                                      x^2 & 0 \\
                                                      0    & x
                                                     \end{pmatrix}
$$
On the other hand we have $(t-a)(t-c)=(t-b)(t-d)$ where
$$
c=\begin{pmatrix}
                 x^4 & 0 \\
                  0 &  x^2
     \end{pmatrix} \quad {\mbox {\rm and } }\quad d=\begin{pmatrix}
                                                      x^4 & x^2 \\
                                                      0    & x^2
                                                     \end{pmatrix}
$$
The ring $A$ is of course right but not left duo (Cf \cite{LTY}
Exercise 22 4A  p.318 ).

}
\end{example}
Let us end this paper with a brief account of some recent
developments related to duo rings and Ore extensions.
\begin{remark}
Hirano, Hong, Kim and Park proved in \cite{HHKP} that an ordinary
polynomial ring is one-sided duo only if it is commutative.  Marks
in \cite{Ma} extended this result to Ore extensions, by showing
that if a noncommutative Ore extension which is  a duo ring on one
side exists,  then it has to be right duo, $\sigma$  must be
non-injective and $\de\ne 0$. He also obtained a series of
necessary conditions for the Ore extension to be right duo.
Matczuk in \cite{M} showed that noncommutative Ore extensions
which are right duo rings do exist and that the necessary
conditions obtained by Marks are not sufficient for the Ore
extension to be right duo.

\end{remark}

\end{document}